\documentclass[xits,review,sort&compress]{elsarticle}
\usepackage[active]{srcltx}
\usepackage{bbm}
\usepackage{mathrsfs}
\usepackage{tipa}
\usepackage{amsthm}
\usepackage{amssymb}
\usepackage{stackrel}
\usepackage{amsmath,amscd}
\usepackage[numbers]{natbib}
\newtheorem{dl}{Theorem}[section]

\newtheorem{yl}[dl]{Lemma}

\newtheorem{tl}[dl]{Corollary}

\newtheorem{xinzhi}[dl]{Proposition}
\newtheorem{remark}[dl]{Remark}

\numberwithin{equation}{section}
\newproof{pot1}{Proof of Theorem \ref{mainthnew}}
\newproof{pot2}{Proof of Theorem \ref{mainthnew2}}
\newproof{pot3}{Proof of Theorem \ref{qanalogue}}
\newproof{pot4}{Proof of Theorem \ref{333}}
\newcommand{\poq}[2]{(#1;q)_{#2}}

\def\pf{\noindent {\it Proof.} }

\begin{document}
\title {Two new transformation formulas for ${}_{8}\psi_{8}$ and ${}_{8}W_{7}$ series  associated with Weierstrass' theta identity}
\author{Jin Wang\fnref{fn3,fn4}}
\fntext[fn3]{This  work was supported by NSF of Zhejiang Province (Grant~No.~LQ20A010004).}
\fntext[fn4]{E-mail address: jinwang@zjnu.edu.cn}
\address[Canada]{Department of Mathematics, Zhejiang Normal  University,
Jinhua 321004,~P.~R.~China}
\author{Xinrong Ma\fnref{fn1,fn2}}
\fntext[fn1]{This  work was supported by NSFC grant No. 11971341.}
\fntext[fn2]{E-mail address: xrma@suda.edu.cn.}
\address[P.R.China]{Department of Mathematics, Soochow University, Suzhou 215006, P. R. China}
\begin{abstract}
In this paper, we establish two new transformation formulas  for ${}_{8}\psi_{8}$ and ${}_8\phi_7$ series by means of  Slater's general transformation for bilateral series. As applications, some specific transformation formulas are presented among which include a general form of Weierstrass' theta identity and new proofs of  Bailey's VWP ${}_6\psi_6$ and Jackson's  ${}_8\phi_7$ summation formula.
\end{abstract}
\begin{keyword} Basic hypergeometric series;  theta function;  Weierstrass' theta identity; transformation; Bailey's ${}_6\psi_6$ summation formula; Jackson's  ${}_8\phi_7$ summation formula
\vspace{5pt}
\end{keyword}
\maketitle
\thispagestyle{empty}
\parskip 7pt

\section{Introduction}
Throughout this paper, we shall adopt the standard notation and terminology for basic
hypergeometric series (or $q$-series) found in Gasper and Rahman's book  \cite{10}. For instance, the $q$-shifted factorial with $0<|q|<1$ and integers $n$ is defined by
\[(x;q)_\infty:
\:=\:\prod_{n=0}^{\infty}(1-xq^n) \quad\text{and}\quad
(x;q)_n:\:=\:\frac{(x;q)_\infty}{(q^nx;q)_\infty}
\]
 and its multi-parameter form by
\[(x_1,x_2,\ldots,x_m;q)_n:\:=\:(x_1;q)_n(x_2;q)_n\cdots(x_m;q)_n.\]

We also need  the  Jacobi theta function
\begin{align}
\theta(x):=\frac{1}{\poq{q}{\infty}}\sum_{n=-\infty}^\infty (-1)^nq^{n(n-1)/2}x^n\quad(x\neq 0)\label{jacobintheta-def}
\end{align}
and the multi-parameter notation
\[\theta(x_1,x_2,\ldots,x_m):\:=\:\theta(x_1)\theta(x_2)\cdots\theta(x_m).\]
Recall that  the classical Jacobi triple product identity (cf. \cite[(II.28)]{10}) asserts
\begin{align}
\theta(x)=(x,q/x;q)_\infty.\label{jacobintheta}
\end{align}
A basic  hypergeometric series with the
base $q$ and the argument $z$   is defined  to be
\begin{align}
{}_{r+1}\phi _{r}\left[\begin{matrix}a_{1},\dots ,a_{r+1}\\ b_{1},\dots ,b_{r}\end{matrix}
; q, z\right]:=\sum _{n=0} ^{\infty }\frac{\poq {a_{1},\ldots
,a_{r+1}}{n}}{\poq
{b_{1},\ldots,b_{r}}{n}}\frac{z^{n}}{\poq{q}{n}},\label{fpfbasic}
\end{align}
while a  bilateral hypergeometric series with the argument $z\neq 0$  is given  by
\begin{eqnarray}
{}_{r}\psi _{r}\left[\begin{matrix}a_{1},\dots ,a_{r}\\
b_{1},\dots ,b_{r}\end{matrix} ; q, z\right]:=\sum _{n=-\infty}
^{\infty }\frac{\poq {a_{1},\ldots,a_{r}}{n}}{\poq{b_{1},\ldots,b_{r}}{n}}z^{n}.\label{fpfbibasic}
\end{eqnarray}
It is convenient for us to employ the shorthand notation
\begin{align}
&\,_{r+1}W_{r}[a_1;a_4,\ldots,a_{r+1};q,z]
\end{align} for the special cases
 of the above ${}_{r+1}\phi _{r}$ series called  \emph{very-well-poised} (in short, VWP), in which all parameters
satisfy the relations
$$
a_1q=a_2b_1=a_3b_2=\cdots=a_{r+1}b_{r}\,\,\mbox{while}\,\, a_2=q\sqrt{a_1},a_3=-q\sqrt{a_1}.
$$
{\bf New notation.}~~In order to state our writing in a convenient form,  we will introduce  two new notation (not used in \cite{10}) as follows:
\begin{align}
\Omega(x_1,x_2,\ldots,x_n;a_1,a_2,\ldots,a_m):=\prod_{i=1}^n\prod_{j=1}^m\poq{x_i/a_j}{\infty}.
\end{align}

Analogously to ${}_{r+1}W_r$ series, we will use the compact notation
\begin{align}{}_{r}\psi _{r}\left[a;a_3,\ldots,a_{r};q,z\right]
\end{align} to denote the  VWP  ${}_{r}\psi_{r}$ series, which means that all parameters are
subject to the relations
$$
a_1b_1=a_2b_2=\cdots=a_{r}b_{r}=aq\,\,\mbox{while}\,\, a_1=q\sqrt{a},a_2=-q\sqrt{a}.
$$

In the context of basic hypergeometric series, finding summation and transformation formulas is an old but interesting study that has attracted
many researchers.   Among the various known results, Weierstrass' theta identity and Bailey's VWP ${}_6\psi_6$ summation formula are very fundamental and frequently used.

\begin{yl}[Weierstrass' theta identity: {\rm cf. \cite{koornwinder}}]\label{weierstrass} For any nonzero complex numbers $x,a,b,c$, there holds
\begin{align}
\theta(xa,x/a,bc,b/c)-\theta(xc,x/c,ab,b/a)=\frac{b}{a}
\theta(xb,x/b,ac,a/c).\label{weierid}
\end{align}
\end{yl}

\begin{yl}[Bailey's VWP  ${}_6\psi_6$ summation formula: {\rm cf. \cite[(II.33)]{10}}]\label{baileypsi66} For any nonzero complex numbers $a,b,c,d,e$ such that $ |a^2q/bcde|<1$, there holds
\begin{align}
   &{}_{6} \psi_{6}[a; b, c, d,e; q, a^2q /bcde]\label{Baileysum-0}\\
&=\frac{(q, q / a,a q, a q / b c, a q / b d, a q / b e, a q / c d, a q / c e, a q / d e; q)_{\infty}}{\left( q / b, q / c, q / d, q / e,a q / b, a q / c, a q / d, a q / e, a^{2}q/ b c d e; q\right)_{\infty}}.\nonumber
\end{align}
\end{yl}

 For the various proofs of Bailey's VWP  ${}_6\psi_6$ summation formula, we refer the reader to R. Askey \cite{andrews}, R. Askey and E. M. Ismail \cite{askey},   W. C. Chu \cite{chu}, F. Jouhet and M. Schlosser \cite{2007},  M. Schlosser \cite{schlosser0,schlosser}, as well as  L. J. Slater and L. Lakin \cite{slater0}. Here we would like to mention that  our previous paper \cite{wangjinphd} reveals a relation of Bailey's VWP  ${}_6\psi_6$ summation formula and Weierstrass' theta identity. Regarding the latter, the reader may consult T. H. Koornwinder's paper \cite{koornwinder} for its history and applications in
the theory of theta functions. Some interesting applications of this identity can be found in the latest paper \cite{wangjinphd-0} by the first author of this paper.

 The theme of the present paper is to establish two new transformations concerning VWP ${}_8\psi_8$
series and  ${}_8W_7$ series. The word ``new" means that these two transformations seem to have not been known in the literature. For instance, they are essentially different from Jackson's transformation formula expressing a VWP ${}_8\psi_8$ series in term of sum of two ${}_8W_7$ series. See \cite[Eq. (5.6.2)]{10} or \cite{jackson} for more details.   For this purpose, we will use  the following result which first appeared in \cite{wei} due to C. A. Wei.

\begin{yl}[Cf. {\rm \cite[Thm. 4]{wei}}]\label{wei-identity} For complex numbers $a,b,c,d,e,f,g$ and $t=bcde/aq$, such that  $|aq/fg|<1, |a^2q/tfg|<1$, there holds
\begin{equation}\label{wei}\begin{aligned}
&{}_{8} \psi_{8}\left[a ; b, c, d, e, f, g ; q, a^{2}q /t fg\right] \\
&=\kappa_0~{}_{8} \psi_{8}[t ; b, c, d, e, t f / a, t g / a ; q, aq / f g]
\\
&+\kappa_1~{}_{8}W_{7}\left[b^{2} / a ; b c / a, b d / a, b e / a, b f / a, b g / a ; q, a^{2}q /t f g\right],
\end{aligned}
\end{equation}
where the coefficients
\begin{subequations}\label{weiid-coeff}
\begin{align}
\kappa_0&:=\frac{\theta(b/a)}{\theta(b/t)}
\frac{\Omega(tq;b,c,d,e)}{\Omega(aq;b,c,d,e)}\frac{\Omega(aq;cd,ce,de,fg,tf,tg)(aq, q / a; q)_{\infty}}{\Omega(b/a;1/c,1/d,1/e)\left(q / f, q / g, tq, q / t, a^{2}q /tfg ; q\right)_{\infty}}, \label{weiid-coeffa}\\
\kappa_1&:=\frac{\theta\left(t/a\right)}{\theta\left(t/b\right)}
\frac{\Omega(bq;b,c,d,e,f,g)}{\Omega(aq;b,c,d,e,f,g)}
\frac{\Omega(b/a;b/ac,b/ad,b/ae)(aq / b f, aq / b g, q / a, aq; q)_{\infty}}{\Omega(b/a;1/c,1/d,1/e)\left(q / f, q / g, q / b, b^{2}q / a ; q\right)_{\infty}}.
\end{align}
\end{subequations}
\end{yl}

Furthermore, we will need to set up another preliminary transformation embedded  in Slater's general transformation for  bilateral series.

\begin{yl}\label{maincoroll-wei}  For complex numbers $a,b,c,d,e,h,k,z$  such that  $|a^3q^2/bcdehk|<1$, there holds
\begin{align}
&\,{}_{8}W_{7}\big[b^2/a;b c  /a,b d  /a,b e  /a,b h  /a,b k  /a;q,a^3 q^2/bcdehk\big]\nonumber\\
&=\lambda_1~{}_{8}\psi_{8}\big[a;b,c,d,e,h,k;q,a^3 q^2/bcdehk\big]\label{wanted}\\
&+\mu_1~{}_{8}\psi_{8}\big[z^2/a;b z/a, c  z/a,d  z/a,e  z/a,h z/a,k z/a;q,a^3 q^2/bcdehk\big],\nonumber
\end{align}
where
\begin{subequations}\label{lambdamu}
\begin{align}\label{lambdamu-1}
 \lambda_1 &=\frac{\poq{b^2 q/a}{\infty}}{\poq{q,aq,q/a}{\infty}}\frac{\Omega(q,aq;b,c,d,e,h,k)}{\Omega(aq/b,bq;c,d,e,h,k)}
\frac{\theta(bz/a,z/b)}{\theta(z/a,z)}, \\
 \mu_1 &=\frac{\poq{b^2 q/a}{\infty}}{\poq{q,z^2 q/a,aq/z^2}{\infty}}
\frac{\Omega(aq/z,zq;b,c,d,e,h,k)}{\Omega(aq/b,bq;c,d,e,h,k)} \frac{\theta(a/b,b)}{\theta(a/z,z)}.
\end{align}
\end{subequations}
\end{yl}

By means of  these two lemmas and Weierstrass' theta identity, we will establish  the following  two new transformation formulas which can be regarded as  generalizations of  Bailey's VWP ${}_6\psi_6$ summation formula and of Weierstrass' theta identity. They are the main results of this paper.

\begin{dl}\label{maintrans-I} For complex numbers $a,b,c,d,e,f,g$ and $z$ with $z=bcde/aq$, such that $|a^2q/fgz|<1, |aq/fg|<1$, there holds
\begin{align}
 &{}_{8} \psi_{8}\left[a ; b, c, d, e, f, g ; q, a^2q/fgz\right] \nonumber\\
 &=(1-1/z)\frac{(aq, q / a; q)_{\infty}}{\Omega(q;f,g)}\frac{\Omega(aq;bc,bd,be,cd,ce,de,fg, fz,gz)}{\Omega(aq;b,c,d,e)\Omega(1/z;1/b,1/c,1/d,
 1/e,fg/a^{2}q )}\nonumber\\
  &\times {}_{8} \psi_{8}[z ; b, c, d, e, fz/ a, gz / a ; q, aq / fg] \label{finalidentities-former}\\
 &+\frac{1}{z}\frac{(aq, q / a; q)_{\infty}}{\left(aq/z^2,z^2q/a;q\right)_{\infty}} \frac{\Omega(aq/z,zq;f,g)}{\Omega(q,aq;f,g)}
 \frac{(b, c, d, e; q)_{\infty}\Omega(aq/z;b,c,d,e)}{\Omega(1/z;1/b,1/c,1/d,1/e)\Omega(aq; b,  c, d,  e)}\nonumber\\
 &\times {}_{8}\psi_{8}\big[z^2/a;bz/a,cz/a,dz/a,
ez/a,fz/a,gz/a;q,a^2q/fgz\big]. \nonumber
\end{align}
 \end{dl}

\begin{dl}\label{maintrans-II}  With the same conditions as Theorem \ref{maintrans-I}. Then there holds
\begin{align}
 &{}_{8}W_{7}\left[b^{2} / a ; b c / a, b d / a, b e / a, b f / a, b g / a ; q, a^{2}q / fgz\right]=(1-1/z)
\frac{\theta(bz/a,z/b)}{\theta(z/a,z)} \nonumber\\
&\times\frac{(b^2 q/a; q)_{\infty}}{\Omega(aq/b, bq;f,g)}\frac{\Omega(q;b,c,d,e)\Omega(aq;f,g,cd,ce,de,fg, fz,gz)}{\Omega( bq;b,c,d,e)\Omega(1/z;1/b,1/c,1/d,
 1/e,fg/a^{2}q )}\nonumber\\
 &\qquad\quad\times {}_{8} \psi_{8}[z ; b, c, d, e, fz/ a, gz / a ; q, aq / f g]\label{finalidentities}\\
 &+\frac{(zq/b;q)_{\infty}}{(z^2q/a;q)_{\infty}}\frac{\theta(1/b)}{\theta(z/a)}
 \frac{\Omega(b/a;1/bq,1/c,1/d,1/e)
 \Omega(zq;f,g)}{\Omega(1/z;z/aq,1/c,1/d,1/e)\Omega(aq/b;f,g)}
\frac{\Omega(aq/z;b,c,d,e,f,g)}{\Omega(bq;b,c,d,e,f,g)}\nonumber\\
 &\qquad\quad\times {}_{8}\psi_{8}\big[z^2/a;bz/a,cz/a,dz/a,
ez/a,fz/a,gz/a;q,a^2q/fgz\big] \nonumber.
\end{align}
 \end{dl}
Our paper is planned as follows. In the next section, we will  show how to derive Lemma \ref{maincoroll-wei} from Slater's
general  transformation which seems to have not previously known. Section 3 is devoted to the proofs of Theorems \ref{maintrans-I} and \ref{maintrans-II}. Their relations with some well-known results such as Bailey's VWP ${}_6\psi_6$ summation formula, Jackson's  ${}_8\phi_7$ summation formula, and Weierstrass' theta identity, are also presented in  Section 4.
\section{The proof of Lemma \ref{maincoroll-wei}}  Lemma \ref{maincoroll-wei} comes from the special case $r=5$ of  Slater's   general transformation (cf. \cite{slater} or \cite[Eq. (5.5.2)]{10}) for VWP ${}_{2r}\psi_{2r}$ series. To make this clear, we first  need to establish
\begin{yl}[Slater's transformation with $r=5$]\label{slatermaintheme}For $x,z\neq  0$ and
$|a^4q^3/bcdefghk|<1$,  it holds
\begin{align}&d_0~{}_{10}W_{9}\big[g^2/a;bg/a,cg/a,dg/a,eg/a,fg/a,
gh/a,gk/a;q,a^4 q^3/bcdefghk\big]\nonumber\\
&=d_1\,{}_{10}\psi_{10}[a;b,c,d,e,f,g,h,k;q,a^4 q^3/bcdefghk]\label{yyy}\\
&+d_2\,{}_{10}\psi_{10}\big[b^2 x^2/a;b^2 x/a,b c x/a,b d x/a,b e x/a,b f x/a,b g x/a,b h x/a,b k x/a;q,a^4 q^3/bcdefghk\big]\nonumber\\
&+d_3~{}_{10}\psi_{10}\big[x^2z^2/a;bx z/a,cx z/a,dx z/a,ex z/a,fx z/a,gx z/a,
hx z/a,kx z/a;q,a^4 q^3/bcdefghk\big],\nonumber
\end{align}
where the coefficients
\begin{subequations}\label{wwwww}
\begin{align}
d_0&=\frac{1}{(g^2 q/a,aq/g^2;q)_{\infty}}\frac{\Omega(gq,aq/g;b,c,d,e,f,g,h,k)}
{\theta(g,b g x/a,g x z/a,a/g,b x/g,x z/g)},\label{16a}\\
  d_1 &= \frac{1}{(a q,q/a;q)_{\infty}}\frac{\Omega(q,aq;b,c,d,e,f,g,h,k)}
  {\theta(g,a/g,b x/a,x z/a,b x,x z)},\\
  d_2 &=\frac{1}{\big(b^2x^2q/a,a q/b^2x^2;q\big)_{\infty}}
  \frac{\Omega(aq/bx,bqx;b,c,d,e,f,g,h,k)}{\theta(b g x/a, b x/g, a/bx, z/b, b x, b x^2 z/a)},\\
  d_3 &=\frac{1}{(x^2 z^2q/a,a q/x^2 z^2;q)_{\infty}}\frac{\Omega(xzq,aq/xz;b,c,d,e,f,g,h,k)}{
  \theta(g x z/a,x z/g,b/z,a/x z,b x^2 z/a,x z)}.
\end{align}
\end{subequations}
\end{yl}
\pf It suffices to set $r=5$ and $b_3=a$ in Slater's transformation   \cite[Eq. (5.5.2)]{10}. Then it is easy to compute that
\begin{align}
d_0\,&{}_{10}W_{9}\left[a; b_{4}, \ldots, b_{10}; q, (aq)^{3}/b_{4}\cdots b_{10}\right]\nonumber\\
&=d_1\, {}_{10}\psi_{10}\left[a_3^2/a; a_{3}, a_{3} b_{4} / a,\ldots, a_{3} b_{10} /a; q, (aq)^{3}/b_{4}\cdots b_{10}\right]\nonumber\\
&+d_2\,
{}_{10}\psi_{10}\left[a_4^2/a;a_{4},  a_{4} b_{4} / a,\ldots, a_{4} b_{10} / a; q, (aq)^{3}/b_{4}\cdots b_{10}\right]\nonumber\\
&+d_3\, {}_{10}\psi_{10}\left[a_5^2/a; a_{5}, a_{5} b_{4} / a, \ldots,a_{5} b_{10} /a; q, (aq)^{3}/b_{4}\cdots b_{10}\right]\label{zzz}
\end{align}
with the corresponding coefficients
\begin{align*}
d_0&=\frac{\left(q,q / b_{4}, \ldots, q / b_{10}, a q / b_{4}, \ldots, a q / b_{10} ; q\right)_{\infty}}{\left(a q; q\right)_{\infty}\theta(a_{3}, a_4, a_{5},a_{3} / a, a_{4} / a, a_{5} / a)},\\
  d_1 &= \frac{\left( q / a_{3},a_{3} q / a,  a_{3} q / b_{4},\ldots, a_{3} q / b_{10},a q / a_{3} b_{4}, \ldots, a q / a_{3} b_{10} ; q\right)_{\infty}}{\left(q a_{3}^{2} / a, a q / a_{3}^{2} ; q\right)_{\infty}\theta(a_{3},a_{3} / a, a_{4} / a_{3},  a_{5} / a_{3},a_{3} a_{4} / a, a_{3} a_{5} / a)},\\
  d_2 &=\frac{\left(q / a_{4},a_{4} q / a, a_{4} q / b_{4}, \ldots, a_{4} q / b_{10}, a q / a_{4}b_{4}, \ldots, a q / a_{4} b_{10} ; q\right)_{\infty}}{\left(q a_{4}^{2} / a, a q / a_{4}^{2} ; q\right)_{\infty}\theta(a_{4},a_{4} / a,  a_{3} / a_{4},  a_{5} / a_{4},a_{4} a_{3} / a, a_{4} a_{5} / a)},\\
  d_3 &=\frac{\left(q / a_{5},a_{5} q / a,a_{5} q / b_{4}, \ldots, a_{5} q / b_{10}, a q / a_{5} b_{4}, \ldots,{a q / a_{5} b_{10}} ; q\right)_{\infty}}{\left(q a_{5}^{2} / a, a q / a_{5}^{2} ; q\right)_{\infty}\theta(a_{5},a_{5} / a, a_{4} / a_{5},  a_{3}/ a_{5},a_{5} a_{4} / a, a_{5}a_{3} / a)}.
\end{align*}
By making the simultaneous substitution
\begin{align*}
&(a,b_4,b_5,b_6,b_7,b_8,b_9,b_{10},a_3,a_4,a_5)\\
&\qquad\qquad\quad\to\big(g^2/a, b g/a,c g/a,d g/a,e g/a, f g/a,g h/a,g k/a, g,b g x/a,g x z/a\big)
\end{align*}
in \eqref{zzz}, then we obtain  \eqref{yyy}
 with the corresponding coefficients given by \eqref{wwwww}. All verifications can be carried by \emph{Mathematica} and left to the interested reader.
\qed

Once specializing  Lemma \ref{slatermaintheme} to  the case $aq=fg$, we come to
\begin{xinzhi}[Contiguous relation for
three VWP ${}_{8}\psi _{8}$ series]\label{mainthtl} With the same notation as Lemma \ref{slatermaintheme}. Then we have
\begin{align}
&\frac{\Omega(q,aq;b,c,d,e,h,k)}{\theta(b x/a,x z/a,b x,x z)}\times{}_{8}\psi_{8}[a;b,c,d,e,h,k;q,a^3 q^2/bcdehk]\label{mainthtl-id}\\
&=\frac{(a q,q/a;q)_{\infty}}{\big(b^2x^2q/a,a q/b^2x^2;q\big)_{\infty}}
  \frac{\Omega(aq/bx,bqx;b,c,d,e,h,k)}{\theta(z/b, b x, bx/a, b x^2 z/a)}\nonumber\\
  &\qquad\qquad\times\,\,{}_{8}\psi_{8}\big[b^2 x^2/a;b^2 x/a,b c x/a,b d x/a,b e x/a,b h x/a,b k x/a;q,a^3 q^2/bcdehk\big]\nonumber\\
&+\frac{(a q,q/a;q)_{\infty}}{(x^2 z^2q/a,a q/x^2 z^2;q)_{\infty}}\frac{\Omega(xzq,aq/xz;b,c,d,e,h,k)}{
  \theta(b/z,x z/a,b x^2 z/a,x z)}\nonumber\\
  &\qquad\qquad\times{}_{8}\psi_{8}\big[x^2 z^2/a;bxz/a, c xz/a,d xz/a,e xz/a,hxz/a,kxz/a;q,a^3 q^2/bcdehk\big].\nonumber
\end{align}
\end{xinzhi}
\pf Note that $d_0=0$ in \eqref{16a} when $aq=fg$.
In this case,  all coefficients given by \eqref{wwwww} can be further simplified. In details, we have
\begin{align*}
  d_1 &= \frac{1}{(a q,q/a;q)_{\infty}}\frac{\Omega(q,aq;b,c,d,e,f,g,h,k)}
  {\theta(g,a/g,bx/a,x z/a,b x,x z)}\\
   &=\frac{\Omega(q,aq;f,g)}{\theta(g,a/g)(a q,q/a;q)_{\infty}}\frac{\Omega(q,aq;b,c,d,e,h,k)}
  {\theta(b x/a,x z/a,b x,x z)}\\
&=\frac{-g/a}{(a q,q/a;q)_{\infty}}\frac{\Omega(q,aq;b,c,d,e,h,k)}
  {\theta(b x/a,x z/a,b x,x z)}.
\end{align*}
By carrying out the same computation, we obtain
\begin{align*}
   d_2 &=\frac{1}{\big(b^2x^2q/a,a q/b^2x^2;q\big)_{\infty}}
  \frac{\Omega(aq/bx,bqx;b,c,d,e,f,g,h,k)}{\theta(b g x/a, b x/g, a/b x, z/b, b x, b x^2 z/a)},\\
  &=\frac{-g/bx}{\big(b^2x^2q/a,a q/b^2x^2;q\big)_{\infty}}
  \frac{\Omega(aq/bx,bqx;b,c,d,e,h,k)}{\theta(a/bx,z/b, b x, b x^2 z/a)};\\
   d_3 &=\frac{1}{(x^2 z^2q/a,a q/x^2z^2;q)_{\infty}}\frac{\Omega(xzq,aq/xz;b,c,d,e,f,g,h,k)}{
  \theta(g x z/a,x z/g,b/z,a/x z,b x^2 z/a,x z)}\\
  &=\frac{-g/xz}{(x^2 z^2q/a,a q/x^2 z^2;q)_{\infty}}\frac{\Omega(xzq,aq/xz;b,c,d,e,h,k)}{
  \theta(b/z,a/x z,b x^2 z/a,x z)}.
\end{align*}
By substituting all these into  \eqref{yyy} and  making   a bit simplification, we obtain  \eqref{mainthtl-id}.
\qed

Having prepared  Proposition \ref{mainthtl}, we are ready to prove  Lemma  \ref{maincoroll-wei}.

\begin{proof}[The proof of Lemma  \ref{maincoroll-wei}]  It suffices to take $x=1$ in  \eqref{mainthtl-id}.
We thus obtain
\begin{align*}
&\frac{\Omega(q,aq;b,c,d,e,h,k)}
  {\theta(b /a,z/a,b,z)}\,\,{}_{8}\psi_{8}[a;b,c,d,e,h,k;q,a^3 q^2/bcdehk]\\
&=\frac{(a q,q/a;q)_{\infty}}{\big(b^2q/a,a q/b^2;q\big)_{\infty}}
  \frac{\Omega(aq/b,bq;b,c,d,e,h,k)}{\theta(z/b, b, b /a, bz/a)}\\
 &\qquad \times{}_{8}\psi_{8}\big[b^2 /a;b^2 /a,b c  /a,b d  /a,b e  /a,b h  /a,b k  /a;q,a^3 q^2/bcdehk\big]\nonumber\\
&+\frac{(a q,q/a;q)_{\infty}}{(z^2q/a,a q/z^2;q)_{\infty}}\frac{\Omega(zq,aq/z;b,c,d,e,h,k)}{\theta(b/z,  z/a,b z/a,  z)}\\
&\qquad\times{}_{8}\psi_{8}\big[z^2/a;b z/a, c  z/a,d  z/a,e  z/a,h z/a,k z/a;q,a^3 q^2/bcdehk\big].\nonumber
\end{align*}
Observe that in such case, the series
$$
{}_{8}\psi_{8}\big[b^2 /a;b^2 /a,b c  /a,b d  /a,b e  /a,b h  /a,b k  /a;q,a^3 q^2/bcdehk\big]
$$
reduces to
$$
{}_{8}W_{7}\big[b^2 /a;b c  /a,b d  /a,b e  /a,b h  /a,b k  /a;q,a^3 q^2/bcdehk\big].
$$
Hence, we can express  the ${}_{8}W_{7}$ series in terms of two VWP ${}_{8}\psi_{8}$ series via the use of the preceding identity. The computation  is carried out as follows:
\begin{align*}
&\frac{(a q,q/a;q)_{\infty}}{\big(b^2q/a,a q/b^2;q\big)_{\infty}}
  \frac{\Omega(aq/b,bq ;b,c,d,e,h,k)}{\theta(z/b, b, b /a, bz/a)}\\
  &\times{}_{8}W_{7}\big[b^2 /a;b c  /a,b d  /a,
  b e /a,b h  /a,b k  /a;q,a^3 q^2/bcdehk\big]\nonumber\\
&=\frac{\Omega(q,aq;b,c,d,e,h,k)}
  {\theta(b /a, z/a,b,z)}\times\,{}_{8}\psi_{8}[a;b,c,d,e,h,k;q,a^3 q^2/bcdehk]\\
&-\frac{(a q,q/a;q)_{\infty}}{(z^2q/a,a q/z^2;q)_{\infty}}\frac{\Omega(zq,aq/z;b,c,d,e,h,k)}{\theta(b/z, z/a,b z/a,  z)}\\
&\times{}_{8}\psi_{8}\big[z^2/a;b z/a, c  z/a,d  z/a,e  z/a,h z/a,k z/a;q,a^3 q^2/bcdehk\big].\nonumber
\end{align*}
It is, after  simplified by  the basic relation of
\(\theta(x)=-x\theta(1/x)\), in agreement with \eqref{wanted}.
\end{proof}

\section{Proofs of Theorems \ref{maintrans-I} and \ref{maintrans-II} }
In this section, we proceed to the  proofs of  our main results.

\subsection{The  proof of Theorem \ref{maintrans-I}} The argument for Theorem \ref{maintrans-I} is entirely based on  Lemmas \ref{weierstrass}, \ref{wei-identity}, and \ref{maincoroll-wei}.

\begin{proof}[The proof of Theorem \ref{maintrans-I}] Actually,  by replacing  $(h,k)$ with $(f,g)$ in  \eqref{wanted} and combining with  \eqref{wei}, we  can come to  a system of linear equations in  $X$ and $Y$
\begin{eqnarray}\label{111-111-111}
\left\{\begin{array}{llll}
&{}_{8} \psi_{8}-\kappa_1&{}_{8}W_{7}&=\kappa_0\,X\\
-\lambda_1&{}_{8} \psi_{8}+&{}_{8}W_{7}&=\mu_1\,Y.
\end{array}
\right.
\end{eqnarray}
For clarity of notation, we now employ the abbreviated  notation ${}_{8} \psi_{8}$ and ${}_{8}W_{7}$  respectively for
 \begin{align}
{}_{8} \psi_{8}&:={}_{8} \psi_{8}\left[a ; b, c, d, e, f, g ; q, a^{3}q^{2}  / b c d e f g\right],\\
{}_{8}W_{7}&:={}_{8}W_{7}\left[b^{2} / a ; b c / a, b d / a, b e / a, b f / a, b g / a ; q, a^{3}q^{2}  / b c d e f g\right]
\nonumber
\end{align}
while $X$ and $Y$ stand,  respectively,  for
\begin{align}
X&:={}_{8} \psi_{8}[t ; b, c, d, e, t f / a, t g / a ; q, aq / f g],
\\
Y&:={}_{8}\psi_{8}\big[z^2/a;bz/a,cz/a,dz/a,
ez/a,fz/a,gz/a;q,a^3q^2/bcdefg\big].\nonumber
\end{align}
Here the coefficients $\kappa_0$ and $\kappa_1$ are the same as in  \eqref{weiid-coeff}, $\lambda_1$ and $\mu_1$ are given by \eqref{lambdamu}.
In what follows, we assume  $$z=t,$$ since it is  under this condition  that we can find a closed-form expression for the   determinant of the system of linear equations \eqref{111-111-111}. The detailed process is as follows:
\begin{align*}
1-\kappa_1\lambda_1
&=1-\frac{\theta(c,d,e,b^2cde/a^2q)}{\theta(bc/a,bd/a,be/a,z)}
\\&=\frac{\theta(bc/a,bd/a,be/a,z)-\theta(c,d,e,b^2cde/a^2q)}
{\theta(bc/a,bd/a,be/a,z)}
\\
&=-\frac{\theta(c,d,e,(aq)^2/b^2cde)-\theta(aq/bc,aq/bd,aq/be,z)}
{\theta(bc/a,bd/a,be/a,z)}.
\end{align*}
In the sequel, we are lucky enough to compute  by appealing  the Weierstrass theta function identity \eqref{weierid}  the numerator
$$
\theta(c,d,e,(aq)^2/b^2cde)-\theta(aq/bc,aq/bd,aq/be,z)=
\frac{aq}{bd}\theta\left(a q/b,b c d/a q,b d e/a q,a q/b c e\right).
$$
Thus we get
\begin{align*}
1-\kappa_1\lambda_1=-\frac{aq}{bd}\frac{\theta\left(a q/b,b c d/a q,b d e/a q,a q/b c e\right)}{\theta(bc/a,bd/a,be/a,z)},
\end{align*}
by which we can solve \eqref{111-111-111} for ${}_{8} \psi_{8}$. That is
\begin{align}
 {}_{8} \psi_{8}&=\frac{\kappa_0}{1-\kappa_1\lambda_1} X+\frac{\kappa_1\mu_1}{1-\kappa_1\lambda_1}Y.\label{IDEN-one}
\end{align}
It remains only to simplify the coefficients. As is expected, we evaluate   under the condition $z=t$ that
\begin{align*}
&\frac{\kappa_0}{1-\kappa_1\lambda_1}=\frac{-bd}{aq}\frac{\theta(bc/a,bd/a,be/a,z)}{\theta\left(a q/b,b c d/a q,b d e/a q,a q/b ce\right)}\\
&\times\frac{(aq, q / a, aq / c d, aq / c e, aq / d e, aq / f g, b / a, tq / c, tq / d, tq / e, aq / t f, aq / t g ; q)_{\infty}}{\left(q / f, q / g, aq / c, aq / d, aq / e, b c / a, b d / a, b e / a, b / t, tq, q / t, a^{2}q / t f g ; q\right)_{\infty}}\\
&=\frac{bd(t-1)}{aq}
\frac{\theta\left(c/t,d/t,e/t\right)}{\theta\left(t/c,
d/t,t/e\right)}\\
&\times\frac{(aq, q / a, aq/bc,aq/bd,aq/be, aq / c d, aq / c e, aq / d e, aq / f g, aq / t f, aq / t g ; q)_{\infty}}
{\left(q / f, q / g, aq/b,aq / c, aq / d, aq / e, b / t, c / t, d / t, e/ t, a^{2}q / t f g ; q\right)_{\infty}}\\
&\stackrel{t=z}{==}\frac{z-1}{z}\frac{(aq, q / a, aq/bc,aq/bd,aq/be, aq / c d, aq / c e, aq / d e, aq / f g, aq/ fz, aq/gz; q)_{\infty}}{\left(q / f, q / g, aq/b,aq / c, aq / d, aq / e, b /z, c /z, d /z, e/z, a^{2}q /fgz ; q\right)_{\infty}}.
\end{align*}
Next,  we evaluate  directly
\begin{align*}
&\frac{\kappa_1\mu_1}{1-\kappa_1\lambda_1}=\frac{\theta\left(z/a\right)}
{\theta\left(z/b\right)}\\&~~\times\frac{(q, aq, q / a, c, d, e, bq/ c, bq/ d, bq/ e, bq/ f, bq/ g, aq / b f, aq / b g ; q)_{\infty}}{\left(q / f, q / g, aq / b, aq / c, aq / d, aq / e, aq / f, aq / g, b c / a, b d / a, b e / a, q / b, b^{2}q / a ; q\right)_{\infty}}\\
&~~\times \frac{\poq{b^2 q/a}{\infty}}{\poq{q,z^2q/a,aq/z^2}{\infty}}
\frac{\Omega(aq/z,zq;b,c,d,e,f,g)}{\Omega(aq/b,bq;c,d,e,f,g)} \frac{\theta(a/b,b)}{\theta(a/z,z)}\times\bigg(-\frac{bd}{aq}\bigg)
\frac{\theta(bc/a,bd/a,be/a,z)}{\theta\left(b/a,z/c, z/e,d/z\right)}\\
&~~=\frac{(aq, q / a,b, c, d, e; q)_{\infty}}{\left(q/f, q/g,z^2q/a,aq/z^2; q\right)_{\infty}}\frac{\Omega(aq/z;b,c,d,e,f,g)\Omega(zq;f,g)}{\Omega( aq; b,  c, d,  e, f, g)\left(b/z,c/z,d/z,e/z;q\right)_\infty}\\
&\qquad\qquad\times\bigg(-\frac{dz}{aq}\bigg)\frac{\theta\left(b/z,c/z,d/z,e/z\right)}{\theta\left(z/b,z/c,d/z,z/e\right)}\\
&~~=\frac{1}{z}\frac{(zq/f,zq/g,aq, q / a,b, c, d, e; q)_{\infty}}
{\left(q /f, q/g,aq/z^2,z^2q/a,b/z,c/z,d/z,e/z;q\right)_{\infty}}
\frac{\Omega(aq/z;b,c,d,e,f,g)}{\Omega(aq; b,  c, d,  e, f, g)}.
\end{align*}
Hence the theorem  is  proved.
\end{proof}

\subsection{The  proof of Theorem \ref{maintrans-II}}

With the help of  \eqref{111-111-111}, we can proceed to the proof of Theorem  \ref{maintrans-II}.

\begin{proof}[The  proof of Theorem \ref{maintrans-II}] Performing as before, we can derive from  \eqref{111-111-111}, under the same condition $z=t$, that
\begin{align}
 {}_{8}W_{7}&=\frac{\lambda_1\kappa_0}{1-\kappa_1\lambda_1} X+\frac{\mu_1}{1-\kappa_1\lambda_1}Y,\label{IDEN-two}
\end{align}
where the coefficient
\begin{align*}
&\frac{\lambda_1\kappa_0}{1-\kappa_1\lambda_1}
=\frac{t-1}{t}\\
&\times\frac{(aq, q / a, aq/bc,aq/bd,aq/be, aq / c d, aq / c e, aq / d e, aq / f g, aq / t f, aq / t g ; q)_{\infty}}{\left(q / f, q / g, aq/b,aq / c, aq / d, aq / e, b / t, c / t, d / t, e/ t, a^{2}q / t f g ; q\right)_{\infty}}\\
&\qquad\times \frac{\poq{b^2 q/a}{\infty}}{\poq{q,aq,q/a}{\infty}}\frac{\Omega(q,aq;b,c,d,e,f,g)}
{\Omega(aq/b,bq;c,d,e,f,g)}
\frac{\theta(bz/a,z/b)}{\theta(z/a,z)}\stackrel{t=z}{==}\frac{z-1}{z}\times
\frac{\theta(bz/a,z/b)}{\theta(z/a,z)}\\&\times\frac{(b^2 q/a,q/b,q/c,q/d,q/e,aq/f,aq/g,aq / c d, aq / c e, aq / d e, aq / f g, aq /fz, aq / g z; q)_{\infty}}{\left(q,aq/bf,aq/bg, bq/c,bq/d,bq/e,bq/f,bq/g,b /z, c /z, d /z, e/z, a^{2}q/f gz; q\right)_{\infty}}.
\end{align*}
In the meantime, we have
\begin{align*}
&\frac{\mu_1}{1-\kappa_1\lambda_1}=-\frac{bd}{aq}
\frac{\theta(bc/a,bd/a,be/a,z)}{\theta\left(b/a,z/c, z/e,d/z\right)}\\
&\qquad\qquad\quad\times\frac{\poq{b^2q/a}{\infty}}{\poq{q,z^2q/a,aq/z^2}{\infty}}
\frac{\Omega(aq/z,zq;b,c,d,e,f,g)}{\Omega(aq/b,bq;c,d,e,f,g)} \frac{\theta(a/b,b)}{\theta(a/z,z)}\\
&~~=\frac{\poq{b^2 q/a}{\infty}}{\poq{q,z^2q/a,aq/z^2}{\infty}}
\frac{\Omega(aq/z,zq;b,c,d,e,f,g)}{\Omega(aq/b,bq;c,d,e,f,g)} \frac{\theta(b)}{\theta(a/z)}\times\frac{d}{q}
\frac{\theta(bc/a,bd/a,be/a)}{\theta\left(z/c, z/e,d/z\right)}\\
&~~=\frac{\poq{b^2 q/a}{\infty}\Omega(zq;b,f,g)}{\poq{q,z^2q/a,aq/z^2}{\infty}\Omega(aq/b;f,g)}
\frac{\Omega(aq/z;b,c,d,e,f,g)}{\Omega(bq;c,d,e,f,g)} \frac{\theta(b)}{\theta(a/z)}
\\&\qquad\qquad\quad\times \frac{d}{q}
\frac{\Omega(zq;c,d,e)}{\theta\left(z/c, z/e,d/z\right)}\times
\frac{\theta(bc/a,bd/a,be/a)}{\Omega(aq/b;c,d,e)}.
\end{align*}
Further simplification yields
\begin{align*}
&\frac{\mu_1}{1-\kappa_1\lambda_1}=\frac{\poq{b^2 q/a,bc/a,bd/a,be/a}{\infty}\Omega(zq;b,f,g)}{\poq{q,z^2q/a,aq/z^2,c/z,d/z,e/z}{\infty}\Omega(aq/b;f,g)}
\frac{\Omega(aq/z;b,c,d,e,f,g)}{\Omega(bq;c,d,e,f,g)} \frac{\theta(b)}{\theta(a/z)}\\
&\qquad\qquad\quad\times\frac{d}{q}
\frac{\theta(c/z,d/z,e/z)}{\theta(z/c, z/e,d/z)}\times
\frac{\theta(bc/a,bd/a,be/a)}{\theta(bc/a,bd/a,be/a)}\\
&~~=\frac{\poq{b^2 q/a,bc/a,bd/a,be/a}{\infty}\Omega(zq;b,f,g)}{\poq{z^2q/a,aq/z^2,c/z,d/z,e/z}{\infty}\Omega(aq/b;f,g)}
\frac{\Omega(aq/z;b,c,d,e,f,g)}{\Omega(bq;b,c,d,e,f,g)} \frac{\theta(1/b)}{\theta(z/a)}.
\end{align*}
It gives the complete proof of the theorem.
\end{proof}

\section{Applications}

As a first interesting application, we can deduce an equivalent version of  Bailey's VWP ${}_6\psi_6$ summation formula  from \eqref{finalidentities}  of  Theorem \ref{maintrans-II} directly.

\begin{tl}  For any nonzero complex numbers $a,b,c,d,e,f,g$ with $z=de/q$  and $a=bc$, such that $ |aq/fg|<1$, there holds
\begin{align}
   &{}_{6} \psi_{6}[z; b, c,  fz / a, gz / a ; q, aq / f g]\label{Baileysum}\\
&=\frac{\left(q,de,q^2/de, de/bc,cq/f,cq/g, bq/f,bq/g,(bcq)^{2}/defg ; q\right)_{\infty}}{(q/b,q/c,aq/f,aq/g,de/b,de/c, aq/fg, aq^2/def, aq^2/deg; q)_{\infty}}\nonumber.
\end{align}
\end{tl}
\pf It suffices to take in \eqref{finalidentities} that $a=bc$, i.e., $zq=de$. As is easily seen, $\poq{bc/a}{\infty}=0$ and
\begin{align*}
 {}_{8}W_{7}[b^{2} / a ; b c / a, b d / a, b e / a, b f / a, b g / a ; q, a^{2}q /fgz]=1.
\end{align*}
Also, from the definitions of $\theta(x)$ and ${}_r\psi_r$, it is not difficult to check that
$$
\frac{\theta(z/a,z)}{\theta(bz/a,z/b)}=
\frac{\theta(a/z,1/z)}{\theta(a/bz,b/z)}
$$
and
$$
{}_{8} \psi_{8}[z ; b, c, d, e, fz / a, gz / a ; q, aq / f g] ={}_{6} \psi_{6}[z ; b, c,  fz / a, gz / a ; q, aq / f g].
$$
All reduces \eqref{finalidentities} to
\begin{align*}
   &{}_{6} \psi_{6}[z ; b, c,  fz / a, gz / a ; q, aq / f g]=\frac{1}{1-q/de}\times
\frac{\theta(bcq/de,q/de)}{\theta(cq/de,bq/de)}\\&\times\frac{\left(q,cq/f,cq/g, bq/c,bq/d,bq/e,bq/f,bq/g,bq/de, cq/de, q/ e, q/d, (bcq)^{2} / defg ; q\right)_{\infty}}{(bq/c,q/b,q/c,q/d,q/e,bcq/f,bcq/g,bq/d, bq/e, bcq /de,  bcq/fg, bcq^2/ def, bcq^2/deg ; q)_{\infty}}\nonumber\\
&=\frac{\left(q,de,q^2/de, de/bc,cq/f,cq/g, bq/f,bq/g,(bcq)^{2}/defg ; q\right)_{\infty}}{(q/b,q/c,bcq/f,bcq/g,de/b,de/c, bcq/fg, bcq^2/ def, bcq^2/deg; q)_{\infty}}\nonumber.
\end{align*}
Thus  \eqref{Baileysum} is proved.
\qed

We remark that under the substitution $(z,b,c,f,g)\to (a,b,c,bcd/a,bce/a)$, \eqref{Baileysum} turns out to be \eqref{Baileysum-0}. As demonstrated in the last section, Weierstrass' theta identity \eqref{wei-identity} plays a very crucial role in the establishments of both Theorem \ref{maintrans-I} and Theorem \ref{maintrans-II}.   Especially noteworthy  that Theorem \ref{maintrans-I}   covers this famous identity conversely.
\begin{tl}[Equivalent form of Lemma \ref{weierstrass}]  For any nonzero complex numbers $a,b,c,d,e$ with $a^2q=bcde$, there holds
\begin{align}
\theta(b,c,d,e)-\theta(a,bc/a,bd/a,be/a)=a\theta\left( b/a,c/a,d/a,e/a\right).\label{finalidentities-former-4}
\end{align}
\end{tl}
\pf It suffices to take $z=a$ in \eqref{finalidentities-former} of Theorem \ref{maintrans-I}. As a result, we obtain
 \begin{align*}
 &{}_{8} \psi_{8}\left[a ; b, c, d, e, f, g ; q,aq/ f g\right] \nonumber\\
 &=\left(1-\frac{1}{a}\right)\frac{(aq, q / a, aq/bc,aq/bd,aq/be, aq / c d, aq / c e, aq / d e, aq / f g, q/f, q/g ; q)_{\infty}}{\left(q / f, q / g, aq/b,aq / c, aq / d, aq / e, b / a, c / a, d / a, e/ a, aq/ f g ; q\right)_{\infty}}\nonumber\\
 &\times {}_{8} \psi_{8}[a ; b, c, d, e, f,g ; q, aq / f g]\\
 &+\frac{1}{a}\frac{(b, c, d, e, q/b, q/c, q/d, q/e; q)_{\infty}}{\left( b/a,c/a,d/a,e/a,aq/b,aq/c,aq/d,aq/e;q\right)_{\infty}}
 \times {}_{8}\psi_{8}\left[a ; b, c, d, e, f, g ; q,aq/ f g\right].\nonumber
\end{align*}
On canceling the common factor ${}_{8}\psi_{8}\left[a ; b, c, d, e, f, g ; q,aq/ f g\right]$ from both sides, we get
\begin{align}
 1 =\frac{(aq, 1/ a, aq/bc,aq/bd,aq/be, aq / c d, aq / c e, aq / d e; q)_{\infty}}{\theta\left( b/a,c/a,d/a,e/a\right)}+\frac{1}{a}\frac{\theta(b,c,d,e)}{\theta\left( b/a,c/a,d/a,e/a\right)}
 .
\label{finalidentities-former-3}
\end{align}
For $bcde=a^2q$,  it is easy to check that
$$
(aq, 1/ a, aq/bc,aq/bd,aq/be, aq / c d, aq / c e, aq / d e; q)_{\infty}=\theta(1/a,bc/a,bd/a,be/a).
$$
Then \eqref{finalidentities-former-4} follows from \eqref{finalidentities-former-3} immediately.
\qed

Beyond this, we also find  a general version of Weierstrass' theta identity from  Theorem \ref{maintrans-I}.

\begin{dl}[Generalized Weierstrass theta identity]\label{dual} For $z=bcde/aq$ such that $|aq/z|<1$, there holds
 \begin{align}
 &{}_{8} \psi_{8}\left[a; b, c, d, e, a/z, z ; q, aq/z\right] \nonumber\\
 &=\bigg(1-\frac{1}{z}\bigg)\frac{(q,q,aq, q /a; q)_{\infty}}{\left(zq/a, q/z,aq/z;q\right)_{\infty}}\frac{\Omega(aq;bc,bd,be,cd,ce,de,z^2)}{\Omega(aq;b,c,d,e)\Omega(1/z;1/b,1/c,1/d,
 1/e)}\nonumber\\
  &\qquad\times{}_{8}W_{7}[1/z ; b/z, c/z, d/z, e/z, z/a; q, q] \label{finalidentities-former-2new}\\
 &+\frac{1}{z}\frac{(q,q,aq, q / a,b, c, d, e; q)_{\infty}}{(q/z,zq/a,zq,aq/z; q)_{\infty}}
 \frac{\Omega(aq/z;b,c,d,e)}{\Omega(aq; b,  c, d,  e)
 \Omega(1/z;1/b,1/c,1/d,1/e)}. \nonumber
\end{align}
 \end{dl}
  \pf  To establish \eqref{finalidentities-former-2new}, it suffices to set $f=a/z,g=z$ in  \eqref{finalidentities-former} of Theorem \ref{maintrans-I} and note that
  $$
  {}_{8}\psi_{8}\big[z^2/a; bz/a,cz/a,dz/a,
ez/a,1,z^2/a;q,aq/z\big]=1.
  $$
  That we obtain is
   \begin{align*}
 &{}_{8} \psi_{8}\left[a; b, c, d, e, a/z, z ; q, aq/z\right] \nonumber\\
 &=\bigg(1-\frac{1}{z}\bigg)\frac{(q,q,aq, q /a; q)_{\infty}}{\left(zq/a, q/z,aq/z;q\right)_{\infty}}\frac{\Omega(aq;bc,bd,be,cd,ce,de,z^2)}
 {\Omega(aq;b,c,d,e)\Omega(1/z;1/b,1/c,1/d,
 1/e)}\nonumber\\
  &\times {}_{8} \psi_{8}[z; b, c, d, e, 1, z^2/a ; q, q] \\
 &+\frac{1}{z}\frac{(q,q,aq, q / a,b, c, d, e; q)_{\infty}}{(q/z,zq/a,zq,aq/z; q)_{\infty}}
 \frac{\Omega(aq/z;b,c,d,e)}{\Omega(aq; b,  c, d,  e)
 \Omega(1/z;1/b,1/c,1/d,1/e)}. \nonumber
\end{align*}
On account of the fact
$$
{}_{8} \psi_{8}[z; b, c, d, e, 1, z^2/a ; q, q] ={}_{8}W_{7}[1/z ; b/z, c/z, d/z, e/z, z/a; q, q],
$$
thus \eqref{finalidentities-former-2new} follows.
  \qed

\begin{remark}
As mentioned earlier,  \eqref{finalidentities-former-2new} covers Weierstrass' theta identity as the special case $z=a$.  The explanation goes as follows.
When $z=a$, the condition becomes $bcde=a^2q$. Thus \eqref{finalidentities-former-2new} is specialized  to
 \begin{align*}
 &{}_{8} \psi_{8}\left[a; b, c, d, e, 1, a ; q, q\right] \nonumber\\
 &=\left(1-\frac{1}{a}\right)\frac{(q,q,aq, q /a,q/a; q)_{\infty}}{\left(q, q,q/a;q\right)_{\infty}}\frac{\Omega(aq;bc,bd,be,cd,ce,de)}{\Omega(aq;b,c,d,e)\Omega(1/a;1/b,1/c,1/d,
 1/e)}\nonumber\\
  &\times{}_{8}W_{7}[1/a ; b/a, c/a, d/a, e/a, 1; q, q]\\
 &+\frac{1}{a}\frac{(q,q,aq, q / a,b, c, d, e; q)_{\infty}}{(q/a,q,aq,q; q)_{\infty}}
 \frac{\Omega(q;b,c,d,e)}{\Omega(aq; b,  c, d,  e)
 \Omega(1/a;1/b,1/c,1/d,1/e)}. \nonumber
\end{align*}
A bit simplification yields
\begin{align*}
1 &=\frac{\theta(1/a,aq/bc,aq/bd,aq/be)}{\theta(b/a,c/a,d/a,e/a)}+ \frac{1}{a}\frac{\theta(b, c, d, e)}{
 \theta(b/a,  c/a, d/a,  e/a)},
\end{align*}
which is identified to be  \eqref{finalidentities-former-4}.
\end{remark}

Beside Theorem \ref{dual}, it deserves to consider other  particular choices of $f$ and $g$, in order to find more transformation formulas.
\begin{tl}\label{maintrans-I-lz} For $z=bcde/aq$, $|aq/z|<1$, there holds
\begin{align}
 &\frac{(aq,1/z; q)_{\infty}\Omega(aq;bc,bd,be,cd,ce,de)}{(b, c, d, e; q)_{\infty}\Omega(aq/z;b,c,d,e)}\times{}_{8} W_{7}[z ; b, c, d, e, z / a; q, q ] \nonumber\\
 &+\frac{1}{z}\frac{\poq{q/z,zq,zq/a,aq/z}{\infty}}{\left(q,q,aq/z^2,z^2q/a;q\right)_{\infty}}\times {}_{8}\psi_{8}\big[z^2/a;bz/a,cz/a,dz/a,
ez/a,z/a,z;q,aq/z\big]\nonumber\\
&=\frac{\Omega(aq;b,c,d,e)\Omega(1/z;1/b,1/c,1/d,
 1/e)}{(b, c, d, e; q)_{\infty}\Omega(aq/z;b,c,d,e)}\label{finalidentities-former-222}.
\end{align}
 \end{tl}
\pf   To show \eqref{finalidentities-former-222}, we only need to take $f=1$ and $g=a$ in \eqref{finalidentities-former} of  Theorem \ref{maintrans-I}. Then we get
\begin{align}
 &{}_{8} \psi_{8}\left[a ; b, c, d, e, 1,a ; q, aq/z\right] \nonumber\\
 &=\bigg(1-\frac{1}{z}\bigg)\frac{(aq, q / a; q)_{\infty}}{\Omega(q;1,a)}\frac{\Omega(aq;bc,bd,be,cd,ce,de,a, z,az)}{\Omega(aq;b,c,d,e)\Omega(1/z;1/b,1/c,1/d,
 1/e,1/aq)}\nonumber\\
  &\times {}_{8} \psi_{8}[z ; b, c, d, e, z / a, z  ; q, q ] \label{finalidentities-former-1}\\
 &+\frac{1}{z}\frac{(aq, q / a; q)_{\infty}}{\left(aq/z^2,z^2q/a;q\right)_{\infty}} \frac{\Omega(aq/z,zq;1,a)}{\Omega(q,aq;1,a)}
 \frac{(b, c, d, e; q)_{\infty}\Omega(aq/z;b,c,d,e)}{\Omega(aq; b,  c, d,  e)\Omega(1/z;1/b,1/c,1/d,1/e)}\nonumber\\
 &\times {}_{8}\psi_{8}\big[z^2/a;bz/a,cz/a,dz/a,
ez/a,z/a,z;q,aq/z\big]. \nonumber
\end{align}
By multiplying both sides of \eqref{finalidentities-former-1} with $\Omega(1/z;1/b,1/c,1/d,1/e)\Omega(aq; b,  c, d,  e)$  and then dividing by $(b, c, d, e; q)_{\infty}\Omega(aq/z;b,c,d,e)$ together some simplifications, we finally obtain \eqref{finalidentities-former-222}.\qed

As a good application of Corollary \ref{maintrans-I-lz}, we can give a new proof of Jackson's $q$-analogue of Dougall's $_{7} F_{6}$ sum  \cite[(II. 22)]{10}.
\begin{xinzhi} For integer $n\geq 0$ and complex numbers $a,b,c,d,z$ with $azq^{n+1}=bcd$, there holds
\begin{align}
{}_{8} W_{7}[z ; b, c, d, z/a, q^{-n}; q, q ]
=\frac{\poq{zq,zq/bc,zq/bd,zq/cd}{n}}{\poq{zq/b,zq/c,zq/d,zq/bcd}{n}}.\label{jacksonphi87}
\end{align}
\end{xinzhi}
\pf  Observe that  \eqref{finalidentities-former-222} is equivalent to
\begin{align}
 &{}_{8} W_{7}[z ; b, c, d, e, z / a; q, q ] \nonumber\\
 &=\frac{1}{z}\frac{\poq{q/z,zq,zq/a,aq/z}{\infty}}{\left(q,q,aq/z^2,z^2q/a;q\right)_{\infty}}
\times\frac{(b, c, d, e; q)_{\infty}\Omega(aq/z;b,c,d,e)}{(aq,1/z; q)_{\infty}\Omega(aq;bc,bd,be,cd,ce,de)}\nonumber
\\ &\times {}_{8}\psi_{8}\big[z^2/a;bz/a,cz/a,dz/a,
ez/a,z/a,z;q,aq/z\big]\nonumber\\
&+\frac{\Omega(aq;b,c,d,e)\Omega(1/z;1/b,1/c,1/d,
 1/e)}{(aq,1/z; q)_{\infty}\Omega(aq;bc,bd,be,cd,ce,de)}.\label{finalidentities-former-555}
\end{align}
In this form, we only need to set $e=q^{-n}$ in \eqref{finalidentities-former-555}. Note that
$(e;q)_\infty=0$ and the condition for Corollary \ref{maintrans-I-lz} becomes $azq^{n+1}=bcd$.  Consequently, \eqref{finalidentities-former-555} can be  simplified as follows:
\begin{align*}
&{}_{8} W_{7}[z ; b, c, d, q^{-n}, z / a; q, q ]\\
 &=\frac{\Omega(bcdq^{-n}/z;b,c,d,q^{-n})\Omega(1/z;1/b,1/c,1/d,
 q^{n})}{(bcdq^{-n}/z,1/z; q)_{\infty}\Omega(bcdq^{-n}/z;bc,bd,cd,bq^{-n},cq^{-n},dq^{-n})}\\
 &=\frac{\poq{cdq^{-n}/z,bdq^{-n}/z,bcq^{-n}/z,bcd/z}{\infty}
 \poq{b/z,c/z,d/z,q^{-n}/z}{\infty}}{(bcdq^{-n}/z,1/z; q)_{\infty}\poq{bq^{-n}/z,cq^{-n}/z,dq^{-n}/z,bc/z,bd/z,cd/z}{\infty}}\\
 &=\frac{\poq{zq,zq/cd,zq/bd,zq/bc}{n}}{\poq{zq/b,zq/c,zq/d,zq/bcd}{n}}.
\end{align*}
Hence \eqref{jacksonphi87} is proved.
\qed

Continuing like above, we have
\begin{tl}\label{3new}For  $a=de$ and $bc=zq$, there holds
 \begin{align}
 {}_{8} \psi_{8}\left[a; b, c, d, e, a/z, z ; q, aq/z\right]=\frac{(q,q,aq, q / a,aq/dz,aq/ez; q)_{\infty}\Delta}{(q/z,zq/a; q)_{\infty}
 \Omega(1/z;1/b,1/c,1/d,1/e)\Omega(aq; b,  c, d,  e)},\label{3newid}
\end{align}
where the factor
\begin{align}
\Delta:=\frac{\theta(d,e)(b, c,aq/bz,aq/cz; q)_{\infty}-\theta(z,a/z)(aq/bd,aq/be,aq/cd,aq/ce;q)_\infty}
{z\,(q/d,q/e,zq,aq/z; q)_{\infty}}.\nonumber
\end{align}
 \end{tl}
\pf  To obtain \eqref{3newid},  we only need to put $bc=zq$ (thus, $a=de$) in  \eqref{finalidentities-former-2new} of  Theorem \ref{dual}. In this case, it is evident that
\[{}_{8}W_{7}[1/z ; b/z, c/z, d/z, e/z, z/a; q, q]={}_{6}W_{5}[1/z; d/z, e/z, z/a; q, q],\]
which can be evaluated by  Rogers' summation formula (II. 20) of \cite{10} in closed-form
\begin{align}
{}_{6}W_{5}[1/z ; d/z, e/z, z/a; q, q]=\frac{\poq{q/z,zq/d e,a q/d z,a q/e z}{\infty}}{\poq{q/d,q/e,aq/de,a q/z^2}{\infty}}.
\end{align}
A direct substitution of  this  into \eqref{finalidentities-former-2new} yields
 \begin{align*}
 &{}_{8} \psi_{8}\left[a; b, c, d, e, a/z, z ; q, aq/z\right] \nonumber\\
 &=\bigg(1-\frac{1}{z}\bigg)\frac{(q,q,aq, q /a; q)_{\infty}}{\left(q/z,zq/a, aq/z;q\right)_{\infty}}\frac{\Omega(aq;bc,bd,be,cd,ce,de,z^2)}{\Omega(aq;b,c,d,e)\Omega(1/z;1/b,1/c,1/d,
 1/e)}\nonumber\\
  &\qquad\qquad\times\frac{\poq{q/z,zq/d e,a q/d z,a q/e z}{\infty}}{\poq{a q/d e,q/d,q/e,a q/z^2}{\infty}}\\
 &+\frac{1}{z}\frac{(q,q,aq, q / a,b, c, d, e; q)_{\infty}}{(q/z,zq/a,zq,aq/z; q)_{\infty}}
 \frac{\Omega(aq/z;b,c,d,e)}{\Omega(aq; b,  c, d,  e)
 \Omega(1/z;1/b,1/c,1/d,1/e)}\\
&=\frac{(q,q,aq, q / a,aq/dz,aq/ez; q)_{\infty}\Delta}{(q/z,zq/a; q)_{\infty}
 \Omega(1/z;1/b,1/c,1/d,1/e)\Omega(aq; b,  c, d,  e)},
\end{align*}
where the factor
\begin{align*}\Delta &:=
\frac{(1-1/z)(1-a/z)\poq{q/z,zq/d e}{\infty}\Omega(aq;bd,be,cd,ce)}{\left(q/d,q/e;q\right)_{\infty}}+\frac{\Omega(aq/z;b,c)(b, c, d, e; q)_{\infty}}{z(zq,aq/z; q)_{\infty}}
 \\
&=\frac{\theta(d,e)(b, c,aq/bz,aq/cz; q)_{\infty}-\theta(z,a/z; q)(aq/bd,aq/be,aq/cd,aq/ce;q)_\infty}{z(q/d,q/e,zq,aq/z; q)_{\infty}}.
\end{align*}
Hence \eqref{3newid} is proved.
\qed
\begin{remark}It should be pointed out that Corollary \ref{3new} is  different with Shukla's VWP ${}_8\psi_8$ summation \cite[Eq. (4.1)]{shukla}, namely,
\begin{align*}{}_8\psi_{8}\left[a; b, c, d, e, z, a q^{2} / z ; q, a^{2}/b c d e\right] =\frac{(q, a q, q / a; q)_{\infty}\Omega(aq;bc,bd,be,cd,ce,de)}{\left( a^{2} q / b c d e ; q\right)_{\infty}\Omega(q,aq; b,  c, d,  e)} \\
\times \left(1-\frac{(1-b c / a)(1-b d / a)(1-b e / a)}{(1-b q / z)(1-b z / a q)\left(1-b c d e / a^{2}\right)}\right) \frac{(1-z / b q)(1-b z / a q)}{(1-z / a q)(1-z / q)}, \end{align*}
because  that  in  \eqref{3newid},  the product  of any two upper-parameters does not equal to $aq^2$.
\end{remark}

\begin{tl}\label{prop49} For $zq=cde$ and $|aq/z|<1$, there holds
\begin{align}
 &{}_{8}W_{7}\left[a; c, d, e, a/z, z ; q, aq/z\right]\nonumber\\
 &=\bigg(1-\frac{1}{z}\bigg)\frac{(q,aq, q /a,aq/z^2; q)_{\infty}}{\theta(a/z)\left(q/z;q\right)_{\infty}}\frac{\Omega(q;c,d,e)\Omega(aq;cd,ce,de)}{\Omega(aq;c,d,e,z)\Omega(1/z;1/c,1/d, 1/e)}\nonumber\\
  &\times{}_{8}W_{7}[1/z ; a/z, c/z, d/z, e/z, z/a; q, q]\label{finished}\\
 &+\frac{1}{z}\frac{(q,aq, q / a, a; q)_{\infty}}{\theta(a/z)(zq; q)_{\infty}}
 \frac{\Omega(1;1/c,1/d,1/e)\Omega(aq/z;c,d,e)}{ \Omega(1/z;1/c,1/d,1/e)\Omega(aq; c, d,  e,z)}. \nonumber
\end{align}
  \end{tl}
\pf  It suffices to take $b=a$ (thus $zq=cde$) in  \eqref{finalidentities-former-2new}, yielding
\begin{align*}
 &{}_{8}W_{7}\left[a; c, d, e, a/z, z ; q, aq/z\right] \nonumber\\
 &=\bigg(1-\frac{1}{z}\bigg)\frac{(q,q,aq, q /a; q)_{\infty}}{\left(zq/a, q/z,aq/z;q\right)_{\infty}}\frac{\Omega(aq;ac,ad,ae,cd,ce,de,z^2)}{\Omega(aq;a,c,d,e)\Omega(1/z;1/a,1/c,1/d, 1/e)}\nonumber\\
  &\times{}_{8}W_{7}[1/z ; a/z, c/z, d/z, e/z, z/a; q, q]\\
 &+\frac{1}{z}\frac{(q,q,aq, q / a,a, c, d, e; q)_{\infty}}{(q/z,zq/a,zq,aq/z; q)_{\infty}}
 \frac{\Omega(aq/z;a,c,d,e)}{ \Omega(1/z;1/a,1/c,1/d,1/e)\Omega(aq; a,  c, d,  e)}. \nonumber
\end{align*}
That is we wanted.\qed

It is of interest to see that Corollary \ref{prop49}  covers a special relation for  ${}_{8}W_{7}$ and VWP ${}_{8}\psi_{8}$ series deserving our attention.
\begin{xinzhi}  For complex numbers $c,d,e$ and $z=cde/q,$  such that $|q/z|<1$, we have
\begin{align}
& {}_{8}W_{7}[1/z ; c/z,d/z,e/z,1/z,z; q, q]\label{rrrrr}\\
&\qquad=\frac{(zq;q)_\infty}{(q/z^2;q)_\infty}\frac{\left(q/z;q\right)_{\infty}^3}{(q; q)_{\infty}^3}\lim_{a\to 1}{}_{8}\psi_{8}[a;\sqrt{a},c,d,e,1/z, z; q, q/z].\nonumber
\end{align}
  \end{xinzhi}
\pf It  only needs to take $a=1$ in \eqref{finished}. Note that $zq=cde$ and $(1;q)_\infty=0$.  Then we have
\begin{align}&\theta(1/z){}_{8}W_{7}[1;c,d,e,1/z, z; q, q/z]\nonumber\\
&=\bigg(1-\frac{1}{z}\bigg)\frac{(q,q, q ,q/z^2; q)_{\infty}}{\left(q/z;q\right)_{\infty}}\frac{\Omega(q;c,d,e)\Omega(q;cd,ce,de)}{\Omega(q;c,d,e,z)\Omega(1/z;1/c,1/d, 1/e)}\nonumber\\
  &\times{}_{8}W_{7}[1/z ; 1/z, c/z, d/z, e/z, z; q, q]\nonumber\\
&=\bigg(1-\frac{1}{z}\bigg)\frac{(q/z^2; q)_{\infty}(q; q)_{\infty}^3}{\left(q/z;q\right)_{\infty}^2}\times{}_{8}W_{7}[1/z ; 1/z, c/z, d/z, e/z, z; q, q].\label{finalresul}
\end{align}
Moreover, it is easily seen that
\begin{align}
&{}_{8}W_{7}[1;c,d,e,1/z, z; q, q/z]
\nonumber\\
&\quad=1+\sum_{n=1}^\infty(1+q^n)\frac{\poq{c,d,e,1/z,z}{n}}{\poq{q/c,q/d,q/e,zq,q/z}{n}}(q/z)^n\label{firstid}\\
&\quad=1+\sum^{-1}_{n=-\infty}(1+1/q^n)\frac{\poq{c,d,e,1/z,z}{-n}}{\poq{q/c,q/d,q/e,zq,q/z}{-n}}(z/q)^n\nonumber\\
&\quad=1+\sum^{-1}_{n=-\infty}(1+q^n)\frac{\poq{c,d,e,1/z,z}{n}}{\poq{q/c,q/d,q/e,zq,q/z}{n}}\bigg(\frac{q^5}{(cde)^2}\bigg)^n(z/q^2)^n\nonumber\\
&\quad=1+\sum^{-1}_{n=-\infty} (1+q^n)\frac{\poq{c,d,e,1/z,z}{n}}{\poq{q/c,q/d,q/e,zq,q/z}{n}}(q/z)^n.\label{secondid}
\end{align}
Therefore, by adding \eqref{firstid} and \eqref{secondid}, we have
\begin{align*}
2&{}_{8}W_{7}[1;c,d,e,1/z, z; q, q/z]\\
&\quad=\sum_{n=-\infty}^\infty(1+q^n)\frac{\poq{c,d,e,1/z,z}{n}}{\poq{q/c,q/d,q/e,zq,q/z}{n}}(q/z)^n,
\end{align*}
which amounts to
\begin{align*}{}_{8}W_{7}[1;c,d,e,1/z, z; q, q/z]=\lim_{a\to 1}{}_{8}\psi_{8}[a;\sqrt{a},c,d,e,1/z, z; q, q/z].
\end{align*}
Substituting this into \eqref{finalresul}  we obtain \eqref{rrrrr} at once.
\qed

We end our paper  with  a concrete transformation  arising from Theorem \ref{maintrans-II}.
\begin{tl}For any nonzero complex numbers $b,c,d,e,f,g$ such that $|bc/fg|<1$ and $|(bc)^{2}/defgq|<1,$ there holds
\begin{align}
 \Omega(bq,c;f,g)&\sum_{n=0}^\infty(1-bq^{2n+1}/c)\frac{\poq{dq/c, eq /c, fq /c, gq /c }{n}}{\poq{bq/d, bq /e, bq /f, bq /g }{n}}\bigg(\frac{(bc)^{2}}{defgq}\bigg)^n\nonumber\\
 &=\kappa_2(b)
\frac{\theta(deq/c,de/b)}{\theta(deq/bc,de)}\frac{\Omega(q;b,c)\Omega(bc;f,g,de,fg, def,deg)}{\Omega(1/de;1/b,1/c,fgq/(bc)^{2})}\nonumber\\
 &\qquad\quad\times {}_{8} \psi_{8}[de; b, c, d, e, defq/bc, degq/bc ; q, bc/ f g]\label{finalidentities-222}\\
 &+\frac{\theta(1/b,deq/b)}{\theta(1/d,1/e)}\frac{\Omega(q/c;1/c,1/d,1/e)\Omega(bc;d f,d g,e f,e g,f g)}{\Omega(1/de;1/bdq,
1/beq,bc/(de)^2q,fgq/(bc)^{2})},\nonumber
 \end{align}
where \[\kappa_2(b):=\frac{(1-b/d)(1-b/e)(1-1/de)}{(1-1/d)(1-1/e)\poq{q}{\infty}}.\]
\end{tl}
\pf Obviously,  when $bc=aq$ and $de=z$,  the series
\begin{align*}
 &{}_{8}W_{7}\left[b^{2} / a ; b c / a, b d / a, b e / a, b f / a, b g / a ; q, a^{2}q / fgz\right]
\nonumber\\
&\quad=\sum_{n=0}^\infty\frac{1-bq^{2n+1}/c}{1-bq/c}\frac{\poq{dq/c, eq /c, fq /c, gq /c }{n}}{\poq{bq/d, bq /e, bq /f, bq /g }{n}}((bc)^{2}/ fgzq)^n;\\
 &{}_{8}\psi_{8}\big[z^2/a;bz/a,cz/a,dz/a,
ez/a,fz/a,gz/a;q,a^2q/fgz\big]\nonumber\\
&\quad={}_{6}\psi_{6}\big[z^2/a;dz/a,
ez/a,fz/a,gz/a;q,a^2q/fgz\big]\\
&\quad=\frac{\poq{q,z^2q/a,aq/z^2,aq/de,aq/d f,aq/d g,aq/e f,aq/e g,aq/f g}{\infty}}{\poq{dq,eq,zq/f,zq/g,aq/dz,aq/ez,aq/fz,aq/g z,a^2 q/f gz}{\infty}}.
\end{align*}
The last equality is given by Bailey's VWP ${}_6\psi_6$ formula. All these reduces \eqref{finalidentities} of   Theorem \ref{maintrans-II} to
\begin{align*}
 &\sum_{n=0}^\infty\frac{1-bq^{2n+1}/c}{1-bq/c}\frac{\poq{dq/c, eq /c, fq /c, gq /c }{n}}{\poq{bq/d, bq /e, bq /f, bq /g }{n}}((bc)^{2}/ fgzq)^n\\
 &=(1-1/z)
\frac{\theta(bz/a,z/b)}{\theta(z/a,z)}\frac{(b^2 q/a; q)_{\infty}}{\Omega(aq/b, bq;f,g)}\frac{\Omega(q;b,c,d,e)\Omega(aq;f,g,cd,ce,de,fg, fz,gz)}{\Omega( bq;b,c,d,e)\Omega(1/z;1/b,1/c,1/d,
 1/e,fg/a^{2}q )}\nonumber\\
 &\qquad\quad\times {}_{8} \psi_{8}[z ; b, c, d, e, fz/ a, gz / a ; q, aq / f g]+\mathbf{S},
 \end{align*}
 where
 \begin{align*}
 \mathbf{S}&:=\frac{(zq/b;q)_{\infty}}{(z^2q/a;q)_{\infty}}\frac{\theta(1/b)}{\theta(z/a)}
 \frac{\Omega(b/a;1/bq,1/c,1/d,1/e)
 \Omega(zq;f,g)}{\Omega(1/z;z/aq,1/c,1/d,1/e)\Omega(aq/b;f,g)}
\frac{\Omega(aq/z;b,c,d,e,f,g)}{\Omega(bq;b,c,d,e,f,g)}\\
&\times \frac{\poq{q,z^2q/a,aq/z^2,aq/d
 e,aq/d f,aq/d g,aq/e f,aq/e g,aq/f g}{\infty}}{\poq{dq,eq,zq/f,zq/g,aq/dz,aq/ez,aq/fz,aq/g z,a^2 q/f gz}{\infty}}.
\end{align*}
All remains to simplify $\mathbf{S}$  by the  conditions $bc=aq$ and $z=de$.  As a result, we achieve
\begin{align*}
 \mathbf{S}=\frac{\theta(1/b,deq/b)}{(1-bq/c)\theta(1/d,1/e)} \frac{\Omega(q/c;1/d,1/e)
 }{\Omega(c;f,g)\Omega(bq;d,e,f,g)}\frac{\poq{q,bc/d f,bc/d g,bc/e f,bc/e g,bc/f g}{\infty}}{\poq{deq/bc,(bc)^2/defgq}{\infty}}.\nonumber
\end{align*}
Finally we are led to
\begin{align}
 &\sum_{n=0}^\infty\frac{1-bq^{2n+1}/c}{1-bq/c}\frac{\poq{dq/c, eq /c, fq /c, gq /c }{n}}{\poq{bq/d, bq /e, bq /f, bq /g }{n}}((bc)^{2}/ defgq)^n\nonumber\\
 &=(1-1/de)
\frac{\theta(deq/c,de/b)}{\theta(deq/bc,de)}\times\frac{(bq^2/c; q)_{\infty}}{\Omega(bq,c;f,g)}\frac{\Omega(q;b,c,d,e)\Omega(bc;f,g,cd,ce,de,fg, def,deg)}{\Omega( bq;b,c,d,e)\Omega(1/de;1/b,1/c,1/d,
 1/e,fgq/(bc)^{2})}\nonumber\\
 &\qquad\quad\times {}_{8} \psi_{8}[de; b, c, d, e, defq/bc, degq/bc ; q, bc/ f g]\label{finalidentities-333}\\
 &+\frac{\theta(1/b,deq/b)}{(1-bq/c)\theta(1/d,1/e)} \frac{\Omega(q/c;1/c,1/d,1/e)
 }{\Omega(bq,c;f,g)}\frac{\poq{bc/d f,bc/d g,bc/e f,bc/e g,bc/f g}{\infty}}{\poq{bq/d,bq/e,deq/bc,(bc)^2/defgq}{\infty}}.\nonumber
 \end{align}
By multiplying both sides of \eqref{finalidentities-333} with $(1-bq/c)\Omega(bq,c;f,g)$ and simplifying the resulted, we obtain  \eqref{finalidentities-222}.
\qed

\begin{remark} It is worth mentioning that the special case $b=f$ of \eqref{finalidentities-222} recovers our results in the paper \cite[Thm. 6 and Thm. 7]{wangjinphd}, claiming that Bailey's VWP ${}_6\psi_6$ summation formula \eqref{Baileysum-0}, Rogers' VWP ${}_6W_5$   summation formula \cite[(II.20)]{10}, and Weierstrass' theta identity \eqref{weierid} are closely related to each others.
\end{remark}

 \end{document}